\documentclass[pdflatex,sn-mathphys-num]{sn-jnl}

\usepackage{amsmath,amssymb,amsthm}
\usepackage{xcolor}

\theoremstyle{thmstyleone}
\newtheorem{theorem}{Theorem}[section]
\newtheorem{proposition}[theorem]{Proposition}
\newcommand{\R}{\mathbb R}
\newcommand{\BJ}{\mathrm{BJ}}
\newcommand{\norm}[1]{\lVert #1\rVert}
\newcommand{\abs}[1]{\lvert #1\rvert}
\newcommand{\inner}[2]{\langle #1,#2\rangle}
\begin{document}

\title[Reverse Sobolev inequalities under Birkhoff--James orthogonality]{Sharp reverse Sobolev inequalities under Birkhoff--James orthogonality}
\hypersetup{
  pdftitle={Sharp reverse Sobolev inequalities under Birkhoff--James orthogonality},
  pdfauthor={Zhiyao Fang, Hui Zhang, Qi Liu, Jinyu Xia}
}

\author[1]{\fnm{Qi} \sur{Liu}}
\email{liuq67@aqnu.edu.cn}

\author*[2]{\fnm{Hui} \sur{Zhang}}
\email{zhangaqtc@126.com}

\author[1]{\fnm{Zhiyao} \sur{Fang}}
\email{fzy73@stu.aqnu.edu.cn}

\author[3]{\fnm{Jinyu} \sur{Xia}}
\email{xiajy23@stu.aqnu.edu.cn}

\affil[1]{\orgdiv{School of Mathematics and Statistics},
  \orgname{Anqing Normal University},
  \orgaddress{\city{Anqing}, \postcode{246133}, \state{Anhui}, \country{China}}}

\affil[2]{\orgdiv{School of Physics and Astronomy},
  \orgname{Anqing Normal University},
  \orgaddress{\city{Anqing}, \postcode{246133}, \state{Anhui}, \country{China}}}
  
\affil[3]{\orgdiv{School of Mathematical Sciences},
  \orgname{Chengdu University of Technology},
  \orgaddress{\city{Chengdu},\postcode{610059}, \state{Sichuan},
  	\country{China}}}

\abstract{We use Birkhoff-James orthogonality to establish reverse Sobolev inequalities and obtain several interesting results concerning these inequalities.}


\maketitle

\section{Introduction}\label{sec:introduction}

A priori estimates for diffusion and incompressible flow often combine an
energy identity with a Sobolev inequality. The energy identity controls
the amount of dissipation in time, and the Sobolev inequality bounds a
Lebesgue norm by the spatial gradient. To continue an evolution from a
bound on a lower-order norm, one may need an estimate in the reverse
direction. In three-dimensional Navier--Stokes flow, for instance, a
uniform bound on the velocity gradient in \(L^2\) allows continuation in
\(H^1\), while the kinetic energy directly controls the velocity in
\(L^2\) \cite[Corollary~5.8]{Tao2013}. The reverse estimates in
Section~\ref{sec:pde} give a way to pass from available norm bounds to
such continuation criteria.

The availability of a reverse estimate depends on the class of states
selected by the equation. In weighted heat-flow estimates, testing the
equation with \(au\) gives the pairing
\(\int_\Omega\nabla u\cdot\nabla(au)\). Its vanishing is equivalent to
stationarity of the weighted square norm \(\int_\Omega au^2\)
\cite{MuruganSaloffCoste2017}. For nonlinear diffusion, the same test
gives \(\int_\Omega|\nabla u|^{p-2}\nabla u\cdot\nabla(au)\), the
corresponding Birkhoff--James condition in \(L^p\). In planar flow,
the scalar vorticity satisfies a transport--diffusion equation, and for
the radial vorticity profile in Section~\ref{sec:pde} the transport
term vanishes pointwise. The diffusion then determines the profile
explicitly and helps identify a weight with zero gradient pairing.
These examples motivate studying how cancellation of a mixed energy
pairing constrains the gradient.

The critical Sobolev embedding into a Lebesgue space is finitely
strictly singular \cite[Theorem~3.3]{LangMihula2023}. The condition
\(\nabla u\perp_{\BJ}\nabla(au)\) selects a homogeneous set of
admissible profiles that need not be a linear subspace. This raises
three variational questions. For a given positive weight, when does
the set contain a nonzero profile, what is the smallest reverse
constant on that set, and when is the constant attained?

For the bounded-domain best-constant problem, let \(\Omega\subset\R^n\) be a bounded connected Lipschitz domain with \(n\ge3\). All functions are real. Write \(H=H_0^1(\Omega)\), \(q=2n/(n-2)\), and \(E(u)=\norm{\nabla u}_2^2\). Norms and integrals refer to \(\Omega\) unless another set is indicated. For \(a\in C^2(\overline\Omega)\) strictly positive, define
\(C_{\BJ}(a,2,\Omega)\)\footnote{This is the least \(C\) for which
\(\norm{\nabla u}_{L^2(\Omega)}\le C\norm{u}_{L^{2n/(n-2)}(\Omega)}\)
for all \(u\) in the displayed constraint. In \(L^2(\Omega;\R^n)\),
Birkhoff--James orthogonality is equivalent to vanishing of the
inner product, which gives \(\int_\Omega\nabla u\cdot\nabla(au)=0\).
The constant is defined whenever this constraint admits a nonzero
function.} by
\begin{equation}\label{eq:definitions}
 C_{\BJ}(a,2,\Omega)
 =\sup_{\substack{0\ne u\in H_0^1(\Omega)\\
                   \int_\Omega \nabla u\cdot\nabla(au)=0}}
 \frac{\bigl(\int_\Omega|\nabla u|^2\bigr)^{1/2}}
      {\bigl(\int_\Omega |u|^{2n/(n-2)}\bigr)^{(n-2)/(2n)}} .
\end{equation}

\begin{samepage}
\textbf{Dimension conventions.}
We use \(n\) for the ambient dimension throughout.
In \(C_{\BJ}(a,2,\Omega)\), the index \(2\) is the gradient exponent.
Theorems~\ref{thm:relaxation} and \ref{thm:whole-space-relaxation}
and the exponential-weight results in Section~\ref{sec:exponential}
use the finite critical exponent \(q=2n/(n-2)\), which requires
\(n>2\). The basic estimate with gradient exponent \(p\) applies when
\(1<p<n\), including \(1<p<2\) in dimension two.

The prescribed-profile criterion in
Proposition~\ref{prop:profile-weight} applies in every dimension
\(n\ge1\). Its sufficient condition based on the sign of the profile
holds in dimensions \(1\le n\le4\). Section~\ref{sec:pde} specifies
the dimension for each example and application. The axisymmetric
and full-velocity continuation results concern three-dimensional
flows. The planar vorticity estimates concern two-dimensional
flows and use any finite \(q>2\).
\end{samepage}

In a normed space, \(x\perp_{\BJ}y\) means
\(\norm{x}\le\norm{x+\tau y}\) for every \(\tau\in\R\)
\cite{James1947}. In the proofs we write
\(Q_a(u)=\int_\Omega\nabla u\cdot\nabla(au)\).
For \(1<p<n\), write \(p^*=np/(n-p)\), take \(u\in W^{1,p}_0(\Omega)\),
and impose the same orthogonality in \(L^p(\Omega;\R^n)\).
If \(m=\min_{\overline\Omega}a>0\) and \(A=\max_{\overline\Omega}a\),
the choice \(\tau=-1/A\) in the defining norm inequality gives
\[
 \norm{\nabla u}_p
 \le \left\|\left(1-\frac aA\right)\nabla u-\frac{u}{A}\nabla a\right\|_p
 \le \left(1-\frac mA\right)\norm{\nabla u}_p+\frac{\norm{u\nabla a}_p}{A}.
\]
Consequently \(m\norm{\nabla u}_p\le\norm{\nabla a}_n\norm{u}_{p^*}\)
by H\"older's inequality.
This gives a finite reverse bound for every \(1<p<n\) on the
orthogonality constraint. Section~\ref{sec:general-lp} gives an exact
supremum formula for the best constant, an attainment criterion in terms
of relaxed maximizers, and the sharp growth order for exponential weights
on cylinders.

The exact formula \eqref{eq:relaxation} and the limits
in \eqref{eq:defect} show that a maximizing sequence can lose gradient
energy only where the weight approaches its minimum.
Proposition~\ref{prop:nonattainment} shows
that this loss prevents attainment on a ball for a smooth positive
radial weight. For exponential weights, the expansion
\eqref{eq:expansion} gives the first correction to the optimal value
at the feasibility threshold, and Theorem~\ref{thm:threshold-expansion}
proves uniqueness of the nearby extremal up to scaling.
Theorem~\ref{thm:whole-space-relaxation} establishes the sharp formula
on \(\R^n\). Compactness of the multiplier \(u\mapsto u\nabla a\)
gives strong convergence in the critical Lebesgue space, while
gradient energy may escape to infinity where the weight approaches
its infimum. Section~\ref{sec:pde} gives an explicit radial
feasibility threshold and a periodic weight with infinite best constant.
For a prescribed profile \(f\), Proposition~\ref{prop:profile-weight}
determines whether a positive weight can enforce exact orthogonality and
how large the ratio \(\sup a/\inf a\) must be.
For axisymmetric Navier--Stokes flow with swirl,
Theorem~\ref{thm:axisym-continuation} gives a continuation criterion
in terms of the weighted gradient pairing of \(u^\theta/r\).
Proposition~\ref{prop:axisym-obstruction} constructs a single smooth
finite-energy datum for which the pairing is positive for every
weight in the criterion. Proposition~\ref{prop:velocity-continuation}
gives the corresponding implication for the full velocity under a
uniform reverse estimate with \(2<q<6\). The examples in
Section~\ref{sec:pde} also treat flow without swirl and planar
vorticity, using an estimate valid for each finite \(q>2\) in the
planar case.

\section{The sharp constant and attainment}\label{sec:relaxation}

Put \(b=\tfrac12\Delta a\) and \(m=\min_{\overline\Omega}a\). Integration
by parts, first for compactly supported smooth functions and then
by density in \(H\), gives
\begin{equation}\label{eq:form}
 Q_a(u)=\int_\Omega a|\nabla u|^2-\int_\Omega bu^2.
\end{equation}
\begin{theorem}\label{thm:relaxation}
Let \(n\ge3\) and let \(\Omega\) and \(a\) satisfy
the bounded-domain hypotheses in Section~\ref{sec:introduction}.
Suppose that some nonzero \(u\in H_0^1(\Omega)\) satisfies
\(\nabla u\perp_{\BJ}\nabla(au)\). The best constant satisfies
\begin{equation}\label{eq:relaxation}
 C_{\BJ}(a,2,\Omega)^2
 =\max_{\substack{0\ne u\in H_0^1(\Omega)\\
                   \int_\Omega a|\nabla u|^2\le
                   \frac12\int_\Omega(\Delta a)u^2}}
 \frac{\tfrac12\int_\Omega(\Delta a)u^2
       -\int_\Omega(a-\min_{\overline\Omega}a)|\nabla u|^2}
      {(\min_{\overline\Omega}a)
       \bigl(\int_\Omega|u|^{2n/(n-2)}\bigr)^{(n-2)/n}} .
\end{equation}
The supremum in \eqref{eq:definitions} is attained if and only if
some maximizer in \eqref{eq:relaxation} satisfies
\(\int_\Omega a|\nabla u|^2=\tfrac12\int_\Omega(\Delta a)u^2\).

Moreover, every maximizing sequence \(u_j\) for \eqref{eq:definitions}
with \(\norm{u_j}_q=1\) has a subsequence converging weakly in \(H\)
and strongly in \(L^q\) to a maximizer \(u\) in \eqref{eq:relaxation}.
On this subsequence,
\begin{equation}\label{eq:defect}
 \begin{aligned}
 \int_\Omega|\nabla(u_j-u)|^2
 &\longrightarrow
 \frac{\tfrac12\int_\Omega(\Delta a)u^2-\int_\Omega a|\nabla u|^2}
      {\min_{\overline\Omega}a},\\
 \int_\Omega(a-\min_{\overline\Omega}a)|\nabla(u_j-u)|^2
 &\longrightarrow 0 .
 \end{aligned}
\end{equation}
\end{theorem}

\begin{proof}
Write \(C=C_{\BJ}(a,2,\Omega)\) and \(F(u)=E(u)-Q_a(u)/m\).
By \eqref{eq:form},
\[
 F(u)=\frac1m\int_\Omega bu^2-\frac1m\int_\Omega(a-m)|\nabla u|^2.
\]
On \(Q_a(u)\le0\), \(mE(u)\le\int_\Omega bu^2\); in particular
imposing \(\norm{u}_q\le1\) gives a bounded set in \(H\).
Weak convergence in \(H\) on a subsequence implies strong convergence
in \(L^2\). Thus \(F\) is weakly upper semicontinuous on this set
and \(Q_a\) is weakly lower semicontinuous. The constraint
\(\norm{u}_q\le1\) is weakly closed by the continuous Sobolev embedding.
The relaxed maximum exists and is positive, by the assumed nonzero
feasible function. Homogeneity gives \(\norm{u}_q=1\) at its maximizers.

For \(Q_a(u)=0\), the numerator of the relaxed quotient is \(E(u)\).
For \(Q_a(u)<0\), take balls \(B_j\Subset\Omega\) whose radii tend to
zero and on which \(\sup_{B_j}(a-m)\to0\). Such balls exist even if
the minimum lies on the boundary, by continuity on \(\overline\Omega\).
Choose \(w_j\in C_c^\infty(B_j)\) such that
\(E(w_j)=1\) and \(\norm{w_j}_q\to0\). To do this, multiply a nonzero
smooth cutoff in each ball by \(\sin(k x_1)\) and normalize its gradient
norm. For each fixed cutoff, the gradient norm grows linearly in \(k\),
by the Riemann--Lebesgue lemma for the square of the leading term, and
the unnormalized \(L^q\) norm remains bounded. The shrinking support
gives \(w_j\rightharpoonup0\) in \(H\) by Cauchy--Schwarz and absolute
continuity of the integral. By \eqref{eq:form}, \(Q_a(w_j)\to m\).
For large \(j\), the quadratic equation \(Q_a(u+s_jw_j)=0\)
has a positive root. Its cross term tends to zero by weak convergence,
so \(s_j\to\sqrt{-Q_a(u)/m}\). Consequently
\[
 E(u+s_jw_j)\longrightarrow F(u),\qquad
 \norm{u+s_jw_j}_q\longrightarrow\norm{u}_q .
\]
This proves \eqref{eq:relaxation}. The attainment criterion follows
from equality of the objectives on the original constraint.

Now let \(u_j\) be the stated normalized maximizing sequence, and
extract a weak \(H\) limit \(u\). Then \(Q_a(u)\le0\), and upper
semicontinuity gives
\[
 0<C^2=\lim_j F(u_j)\le F(u)\le C^2\norm{u}_q^2\le C^2.
\]
All inequalities are equalities. Uniform convexity upgrades weak
\(L^q\) convergence with equal limiting norms to strong convergence.
Weak convergence of the gradients gives
\(E(u_j-u)\to C^2-E(u)=-Q_a(u)/m\).
Since \(F(u_j)\to F(u)\) and \(u_j\to u\) in \(L^2\), the weighted
gradient norms for the weight \(a-m\) converge as well. Their cross
term converges by weak convergence, proving the second limit in
\eqref{eq:defect}.
\end{proof}

The second limit in \eqref{eq:defect} locates any missing gradient
energy where \(a\) approaches its minimum. For example, on
\(\{a\ge m+\eta\}\) with fixed \(\eta>0\), the gradients converge
strongly in \(L^2\).

The substitution \(v=\sqrt a\,u\) identifies the weights for which
the orthogonality constraint is feasible. Directly,
\begin{equation}\label{eq:gauge}
 \int_\Omega\nabla u\cdot\nabla(au)
 =\int_\Omega|\nabla v|^2-\frac14\int_\Omega|\nabla\log a|^2 v^2 .
\end{equation}
Indeed, with \(g=\tfrac12\log a\), the gradients of \(u\) and \(au\)
are \(e^{-g}(\nabla v-v\nabla g)\) and \(e^g(\nabla v+v\nabla g)\).
Multiplication by \(e^{\pm g}\) is bounded on \(H\). Nonzero feasible
functions therefore exist exactly when the lowest Dirichlet eigenvalue
of \(-\Delta-\tfrac14|\nabla\log a|^2\) is nonpositive. If it is zero,
all feasible functions are obtained by multiplying the lowest
eigenspace by \(a^{-1/2}\). If it is negative, combine a negative
eigenfunction with one for a positive eigenvalue. A positive
eigenvalue exists because the potential is bounded and the Dirichlet
eigenvalues diverge.

\begin{proposition}\label{prop:nonattainment}
Let \(B=B_1(0)\subset\R^3\) and \(a_m(x)=m+|x|^2\). Set
\begin{equation}\label{eq:D}
 D=\sup_{\substack{v\in H_0^1(B)\\\norm{v}_6=1}}
       \left(3\int_B v^2-\int_B |x|^2|\nabla v|^2\right).
\end{equation}
Then \(D>0\). For all sufficiently small \(m>0\), the constraint is
feasible and \(C_{\BJ}(m+|x|^2,2,B_1(0))=\sqrt{D/m}\), with no extremal.
\end{proposition}

\begin{proof}
Write \(r=|x|\), and denote the expression in parentheses by \(P(v)\).
For the radial function \(f=r^{-4/3}-1\), the one-dimensional
integrals are
\[
 (4\pi)^{-1}\int_B f^2=\frac{32}{15},\qquad
 (4\pi)^{-1}\int_B r^2|\nabla f|^2=\frac{16}{3}.
\]
Cutting \(f\) off in \(r<2\varepsilon\) approximates these two integrals.
Indeed, the gradient of a radial cutoff is bounded by \(C/\varepsilon\)
on \(\varepsilon<r<2\varepsilon\), where its weighted squared
contribution is bounded by \(C\int_{r<2\varepsilon} f^2\to0\).
The cutoff functions belong to \(H_0^1(B)\), and \(P(f)>0\) by the
displayed values. Normalizing a cutoff in \(L^6\) proves \(D>0\).

To prove attainment in \eqref{eq:D}, first maximize \(P\) on
\(\norm{v}_6\le1\) in the space
\[
 \mathcal X=
 \{v\in L^6(B): v\in H^1(B\setminus\overline B_\varepsilon)
 \text{ for every }0<\varepsilon<1,\ 
 v|_{\partial B}=0,\ r\nabla v\in L^2(B)\}.
\]
The cutoff argument just used also approximates every \(v\in\mathcal X\)
in \(L^6\) and in weighted gradient norm. Thus the maximum sought has
the same supremum \(D\). A maximizing sequence has bounded weighted
gradient norms and is bounded in \(H^1\) on each annulus. Diagonal
extraction gives a weak \(L^6\) limit, weak \(H^1\) limits on annuli,
and strong \(L^2\) convergence there. The bound
\(\int_{B_\varepsilon} |v|^2\le |B_\varepsilon|^{2/3}\norm{v}_6^2\)
upgrades the latter to strong \(L^2(B)\) convergence. The zero trace
passes to the limit on an outer annulus. Lower semicontinuity
of the weighted gradient integral, by exhaustion, shows that the
limit maximizes \(P\). Since \(D>0\), its \(L^6\) norm is one.

Let \(v\) be any such maximizer. Replacing it by \(\abs v\) preserves
its norm and gradient energy, so take \(v\ge0\) in the following
estimates. Differentiating \(P(v)/\norm{v}_6^2\) in directions
\(h\in\mathcal X\) gives
\begin{equation}\label{eq:weightedEL}
 \int_B r^2\nabla v\cdot\nabla h=\int_B (3v-Dv^5)h.
\end{equation}
Set \(M=(3/D)^{1/4}\). Testing with \((v-M)_+\in\mathcal X\) gives
a nonnegative left-hand side, while the integrand on the right is
strictly negative wherever \(v>M\). Hence \(v\le M\) almost everywhere.
On a shell \(s<r<2s\) with \(0<s<1/4\), choose \(\eta=1\) on the shell
with support in \(s/2<r<3s\) and \(\norm{\nabla\eta}_\infty\le C/s\).
In \eqref{eq:weightedEL} test with \(\eta^2v\) and absorb the cross
term by Young's inequality. This yields
\[
 \int_B r^2\eta^2|\nabla v|^2
 \le 6\int_B\eta^2v^2+4\int_B r^2|\nabla\eta|^2v^2
 \le C M^2 s^3 .
\]
In particular \(\int_{s<r<2s}|\nabla v|^2\le C M^2s\).
Summing over dyadic shells controls the unweighted gradient near the
origin. Away from the origin, the uniform weighted energy bound
\(\int_B r^2|\nabla v|^2\le 3|B|^{2/3}\) suffices.
Cutoffs at the origin now converge in unweighted \(H^1\), because the
extra unweighted gradient term has squared norm at most
\(C M^2\varepsilon\). This proves \(v\in H_0^1(B)\). The same bound
holds for a signed maximizer: the energy of \(\abs v\) agrees with
that of \(v\) on annuli, and the cutoff argument applies with
\(\abs v\le M\). Therefore all normalized maximizers in
\eqref{eq:D} have \(E(v)\le K\) for a fixed finite \(K>0\).

Now take \(0<m<D/K\). For any normalized maximizer in \eqref{eq:D},
\(Q_{a_m}(v)=mE(v)-D<0\). There is a positive direction for \(Q_{a_m}\)
by \eqref{eq:gauge}, so the equality constraint is nontrivial.
Since \(E(v)-Q_{a_m}(v)/m=P(v)/m\), \eqref{eq:relaxation} gives
\(C_{\BJ}(m+|x|^2,2,B_1(0))^2=D/m\). Every relaxed maximizer has \(Q_{a_m}<0\) after
normalization, so the attainment criterion excludes an extremal.
\end{proof}

On \(\R^n\), let \(\dot W^{1,p}(\R^n)\) be the completion of
\(C_c^\infty(\R^n)\) for the gradient norm. We identify its elements with
their representatives in \(L^{np/(n-p)}(\R^n)\) using the Euclidean
Sobolev inequality \cite{Talenti1976}, and write \(\dot H^1=\dot W^{1,2}\).
Suppose \(a\in C^1(\R^n)\), \(0<m=\inf_{\R^n}a\le a\le M<\infty\),
and \(\nabla a\in L^n(\R^n)\). Multiplication by \(a\) maps
\(\dot W^{1,p}\) into itself by the product rule and H\"older's inequality.
For \(1<p<n\), differentiating the \(p\)-th power of the gradient norm
along \(\nabla(au)\) gives
\[
 \nabla u\perp_{\BJ}\nabla(au)
 \quad\Longleftrightarrow\quad
 \int_{\R^n}a|\nabla u|^p\,dx
   +\int_{\R^n}u|\nabla u|^{p-2}\nabla u\cdot\nabla a\,dx=0.
\]
The second integral is bounded in absolute value by
\(\norm{\nabla u}_{L^p}^{p-1}
  \norm{u}_{L^{np/(n-p)}}\norm{\nabla a}_{L^n}\). In particular,
\begin{equation}\label{eq:whole-space-reverse}
 \displaystyle
 \norm{\nabla u}_{L^p(\R^n)}
 \le \frac{\norm{\nabla a}_{L^n(\R^n)}}{\inf_{\R^n}a}\,
       \norm{u}_{L^{np/(n-p)}(\R^n)}
 \quad\text{if }\nabla u\perp_{\BJ}\nabla(au).
\end{equation}

\begin{theorem}\label{thm:whole-space-relaxation}
Let \(n\ge3\), so that the critical exponent
\(2n/(n-2)\) is finite, and let \(a\) satisfy the preceding whole-space
hypotheses. Suppose there is a nonzero \(u\in\dot H^1(\R^n)\) with
\(\int_{\R^n}\nabla u\cdot\nabla(au)\,dx=0\).
Define \(C_{\BJ}(a,2,\R^n)\) by the supremum of the reverse Sobolev
quotient over these nonzero functions. Then
\begin{equation}\label{eq:whole-space-relaxation}
 \displaystyle
 C_{\BJ}(a,2,\R^n)^2
 =\max_{\substack{0\ne u\in\dot H^1(\R^n)\\
                  \int_{\R^n}\nabla u\cdot\nabla(au)\,dx\le0}}
 \frac{\displaystyle
       \int_{\R^n}|\nabla u|^2\,dx
       -m^{-1}\int_{\R^n}\nabla u\cdot\nabla(au)\,dx}
      {\displaystyle
       \left(\int_{\R^n}|u|^{2n/(n-2)}\,dx\right)^{(n-2)/n}} .
\end{equation}
The orthogonality supremum is attained exactly when some maximizer
on the right satisfies \(\int_{\R^n}\nabla u\cdot\nabla(au)\,dx=0\).
Every maximizing sequence for the orthogonality supremum with
\(\norm{u_j}_{L^{2n/(n-2)}}=1\) has a subsequence converging weakly
in \(\dot H^1\) and strongly in \(L^{2n/(n-2)}\) to a relaxed
maximizer \(u\). Along that subsequence,
\begin{align}
 \int_{\R^n}|\nabla(u_j-u)|^2\,dx
   &\longrightarrow -m^{-1}\int_{\R^n}\nabla u\cdot\nabla(au)\,dx,
       \label{eq:whole-space-gradient-defect}\\
 \int_{\R^n}(a-m)|\nabla(u_j-u)|^2\,dx
   &\longrightarrow 0.
       \label{eq:whole-space-weighted-defect}
\end{align}
\end{theorem}

\begin{proof}
Write \(q=2n/(n-2)\), \(E(u)=\int_{\R^n}|\nabla u|^2\,dx\),
\(Q_a(u)=\int_{\R^n}\nabla u\cdot\nabla(au)\,dx\) and
\(F(u)=E(u)-Q_a(u)/m\). The map \(u\mapsto u\nabla a\) from
\(\dot H^1\) to \(L^2\) is compact. Indeed, a bounded weakly convergent
sequence in \(\dot H^1\) converges strongly in \(L^2\) on each bounded
ball by the Sobolev bound and the local Rellich theorem. Approximating
\(\nabla a\) in \(L^n\) by bounded compactly supported vector fields
gives the claimed compactness, since
\(\norm{u(\nabla a-b)}_2\le\norm{u}_q\norm{\nabla a-b}_n\).
In particular, \(\int u\nabla a\cdot\nabla u\) is weakly continuous.
The identities
\[
 Q_a(u)=\int a|\nabla u|^2+\int u\nabla a\cdot\nabla u,\qquad
 F(u)=-m^{-1}\int(a-m)|\nabla u|^2
      -m^{-1}\int u\nabla a\cdot\nabla u
\]
give weak lower semicontinuity of \(Q_a\) and weak upper
semicontinuity of \(F\). On \(Q_a(u)\le0\), the product estimate
used in \eqref{eq:whole-space-reverse} gives
\(m\norm{\nabla u}_2\le\norm{\nabla a}_n\norm u_q\).
Thus \(F\) has a maximum on \(Q_a(u)\le0,\ \norm u_q\le1\).
It is positive by the nonzero orthogonal profile, so all its
maximizers have \(\norm u_q=1\).

For \(Q_a(u)<0\), choose balls \(B_j\) of radii less than \(1/j\)
on which \(\sup_{B_j}(a-m)<1/j\), using continuity and the definition
of \(m\). On each \(B_j\), a smooth cutoff times a sufficiently
high-frequency sine, normalized in gradient norm, yields
\(w_j\in C_c^\infty(B_j)\) with \(E(w_j)=1\) and \(\norm{w_j}_q\le1/j\).
The support measures tend to zero, hence \(w_j\rightharpoonup0\)
in \(\dot H^1\), even if the centers escape to infinity. The product
formula gives \(Q_a(w_j)\to m\), since the cross term is bounded by
\(\norm{w_j}_q\norm{\nabla a}_n\). The mixed term in
\(Q_a(u+s w_j)\) tends to zero for fixed \(s\) by weak convergence
and continuity of the associated bilinear form. Its positive root
\(s_j\) therefore tends to \(\sqrt{-Q_a(u)/m}\).
Consequently \(E(u+s_jw_j)\to F(u)\) and
\(\norm{u+s_jw_j}_q\to\norm u_q\). On \(Q_a(u)=0\), \(F(u)=E(u)\).
This proves \eqref{eq:whole-space-relaxation} and the attainment criterion.

For a normalized maximizing sequence with \(Q_a(u_j)=0\), extract a
weak \(\dot H^1\) limit \(u\). Lower semicontinuity gives \(Q_a(u)\le0\).
By upper semicontinuity and \eqref{eq:whole-space-relaxation},
\[
 C_{\BJ}(a,2,\R^n)^2
 \le F(u)\le C_{\BJ}(a,2,\R^n)^2\norm u_q^2
 \le C_{\BJ}(a,2,\R^n)^2.
\]
Thus \(u\) is a relaxed maximizer and \(\norm u_q=1\).
Weak \(L^q\) convergence and uniform convexity yield strong \(L^q\)
convergence. The gradient norm satisfies
\(E(u_j-u)\to F(u)-E(u)=-Q_a(u)/m\) by weak convergence.
Finally, convergence of \(F(u_j)\) and weak continuity of the mixed
integral imply convergence of the \((a-m)\)-weighted gradient norms.
Weak convergence of these weighted gradients then gives
\eqref{eq:whole-space-weighted-defect}.
\end{proof}

If \(a>m\) everywhere, \eqref{eq:whole-space-weighted-defect} gives
strong convergence of the gradients on each compact set. In that
case any positive defect in \eqref{eq:whole-space-gradient-defect}
escapes to infinity.

We now prescribe the profile and seek a weight. For any integer \(n\ge1\),
let \(\mathcal A_n\) consist of the smooth bounded functions on \(\R^n\)
with \(\inf a>0\) and compactly supported gradient. Write
\(s_+=\max\{s,0\}\) and \(s_-=\max\{-s,0\}\).

\begin{proposition}\label{prop:profile-weight}
Let \(n\ge1\) and let \(0\ne f\in H^2(\R^n)\) be real. There exists
a weight \(a\in\mathcal A_n\) such that
\(\nabla f\perp_{\BJ}\nabla(af)\) in \(L^2(\R^n;\R^n)\) if and only
if \(\int_{\R^n}(f\Delta f)_+\,dx>0\).
When such a weight exists,
\begin{equation}\label{eq:profile-contrast}
 \displaystyle
 \inf_{\substack{a\in\mathcal A_n\\ \nabla f\perp_{\BJ}\nabla(af)}}
       \frac{\sup_{\R^n}a}{\inf_{\R^n}a}
   =\frac{\int_{\R^n}(f\Delta f)_-\,dx}
          {\int_{\R^n}(f\Delta f)_+\,dx}.
\end{equation}
Every nonzero \(f\in C^2\cap H^2(\R^n)\) of one sign admits such a
weight if \(1\le n\le4\). In every \(n\ge5\) there are smooth positive
\(H^2\) profiles that admit no such weight.
For an axisymmetric profile in \(\R^3\) the weights in the infimum
may also be taken axisymmetric.
\end{proposition}

\begin{proof}
All integrals in this proof are on \(\R^n\) unless
a smaller set is indicated. Put \(h=f\Delta f\in L^1\) and \(E=\int|\nabla f|^2>0\).
Integration by parts gives \(\int h=-E\) and
\(\int\nabla f\cdot\nabla(af)=-\int ah\).
Orthogonality implies
\((\sup a)\int h_+\ge(\inf a)\int h_-\), giving necessity and the
lower bound in \eqref{eq:profile-contrast}.
If \(\int h_+>0\), take \(0\le\chi_j\le1\) in \(C_c^\infty(\R^n)\)
approaching \(\mathbf1_{\{h>0\}}\) in \(L^1(|h|\,dx)\). For large
\(j\), \(\int\chi_jh>0\), and the explicit weight
\begin{equation}\label{eq:profile-cutoff}
 a_j(x)=1+\frac{\int_{\R^n}|\nabla f|^2\,dx}
                     {\int_{\R^n}\chi_j f\Delta f\,dx}\,\chi_j(x)
\end{equation}
satisfies \(\int a_jh=0\). Its ratio is at most
\(1+E/\int\chi_jh\), which tends to \(\int h_-/\int h_+\).
Rotation averaging of the cutoffs preserves \(\int\chi_jh\) for
axisymmetric \(f\).

For the assertion about profiles of one sign, take \(f\ge0\) after
changing its sign if necessary. If \(f\Delta f\le0\), then
\(-\Delta f\ge0\) on the closure of \(\{f>0\}\) by continuity,
and \(\Delta f=0\) in the interior of the zero set. Let \(\zeta\)
be a nonnegative smooth cutoff equal to one on \(B_1\) and zero
outside \(B_2\), and put \(\zeta_R(x)=\zeta(x/R)\).
The support of \(\Delta\zeta_R\) lies in the annulus, so
\[
 0\le\int_{B_R}(-\Delta f)\,dx
 \le\int f(-\Delta\zeta_R)\,dx
 \le R^{n/2-2}\norm{f}_{L^2(B_{2R}\setminus B_R)}
                  \norm{\Delta\zeta}_2 .
\]
For \(1\le n\le4\), the last expression tends to zero by
\(f\in L^2\), including at \(n=4\). Thus \(f\) is an \(L^2\)
harmonic function and hence vanishes, a contradiction.

For \(n\ge5\), the positive function
\(f(x)=(1+|x|^2)^{-(n-2)/2}\) belongs to \(H^2(\R^n)\), since it
decays as \(|x|^{2-n}\) and its first and second derivatives decay
faster. It satisfies
\[
 -\Delta f=n(n-2)(1+|x|^2)^{-2}f,
\]
so \((f\Delta f)_+=0\). The weight criterion excludes this profile.
\end{proof}

\section{Sharp constants for exponential weights}\label{sec:exponential}

Throughout this section \(n\ge3\) and \(q=2n/(n-2)\).
Set \(a_t(x)=e^{2tx_1}\) for \(t\ge0\).
Let \(\lambda_1<\lambda_2\) be the first two Dirichlet
eigenvalues, with the first simple by connectedness, and let
\(\phi>0\) be the \(L^2\)-normalized first eigenfunction. In
\eqref{eq:gauge}, \(v=e^{tx_1}u\) and \(Q_{a_t}(u)=E(v)-t^2\norm{v}_2^2\).
The constraint is feasible exactly for \(t\ge t_0=\sqrt{\lambda_1}\).
At equality, the set of feasible functions is the span of
\(U=e^{-t_0x_1}\phi\). For larger \(t\) the constraint can be realized
by mixing \(\phi\) with any eigenfunction whose eigenvalue exceeds \(t^2\).
The weak equation for \(\phi\), tested with \(e^{-2t_0x_1}\phi\),
gives \(E(U)=(\lambda_1+t_0^2)\norm{U}_2^2\). Hence
\begin{equation}\label{eq:threshold-value}
 C_{\BJ}(e^{2\sqrt{\lambda_1}\,x_1},2,\Omega)
 =\sqrt{2\lambda_1}\,
 \frac{\bigl(\int_\Omega e^{-2\sqrt{\lambda_1}\,x_1}\phi^2\bigr)^{1/2}}
      {\bigl(\int_\Omega
       |e^{-\sqrt{\lambda_1}\,x_1}\phi|^{2n/(n-2)}
       \bigr)^{(n-2)/(2n)}} .
\end{equation}

\begin{theorem}\label{thm:threshold-expansion}
There is a \(\delta>0\) such that for \(0<t^2-\lambda_1<\delta\)
the constant \(C_{\BJ}(e^{2tx_1},2,\Omega)\) is attained by a unique extremal up to
nonzero scalar multiplication. There is a \(K_\Omega>0\) such that
\begin{equation}\label{eq:expansion}
 \begin{aligned}
 C_{\BJ}(e^{2tx_1},2,\Omega)
 &=C_{\BJ}(e^{2\sqrt{\lambda_1}\,x_1},2,\Omega)\\
 &\quad+K_\Omega\sqrt{t^2-\lambda_1}+O(t^2-\lambda_1)
 \qquad (t\downarrow\sqrt{\lambda_1}).
 \end{aligned}
\end{equation}
To describe \(K_\Omega\), put
\[
 X=\{h\in H:\inner{h}{\phi}_{L^2}=0\},\qquad
 B(h,k)=\int_\Omega\nabla h\cdot\nabla k-\lambda_1\int_\Omega hk,
 \qquad M=e^{-t_0x_1}.
\]
Then \(K_\Omega=\norm{g}_B\), where \(g\in X\) is determined by
\begin{equation}\label{eq:riesz}
 B(g,h)=C_{\BJ}(e^{2\sqrt{\lambda_1}\,x_1},2,\Omega)\left(
 \frac{\int_\Omega\nabla U\cdot\nabla(Mh)}{\int_\Omega|\nabla U|^2}
 -\frac{\int_\Omega U^{(n+2)/(n-2)}Mh}{\int_\Omega U^{2n/(n-2)}}\right)
 \quad (h\in X).
\end{equation}
\end{theorem}

\begin{proof}
Let \(s=\sqrt{t^2-\lambda_1}\) and assume \(0<s^2<\lambda_2-\lambda_1\).
By the spectral gap \(B\) is an inner product on \(X\) whose norm
is equivalent to \(\norm{\nabla\cdot}_2\). Every nonzero feasible
\(v=e^{tx_1}u\) has nonzero projection onto \(\phi\), since the Rayleigh
quotient on \(X\setminus\{0\}\) is at least \(\lambda_2\). Divide by its
\(\phi\) coefficient to write \(v=\phi+h\), \(h\in X\). The constraint
becomes \(B(h,h)=s^2(1+\norm{h}_2^2)\), so it is equivalent to
\begin{equation}\label{eq:parameterization}
 h=h_s(z):=\frac{sz}{\sqrt{1-s^2\norm{z}_2^2}},\qquad \norm{z}_B=1.
\end{equation}
The inverse map is \(z=h/\norm{h}_B\).

Let \(R(u)=\sqrt{E(u)}/\norm{u}_q\) on \(H\setminus\{0\}\).
The expression on the right of \eqref{eq:riesz} is
\(\ell(h)=DR(U)[Mh]\), a bounded functional on \(X\). To prove that
\(\ell\ne0\), suppose \(\ell=0\). Since \(R\) is homogeneous of degree
zero along positive scalings, \(DR(U)[U]=0\) as well. The decomposition
\(H=\R\phi\oplus X\) and the invertibility of multiplication by \(M\)
give \(DR(U)=0\) on \(H\). Hence
\[
 -\Delta U=cU^{q-1}\quad\text{in distributions},\qquad
 c=E(U)/\norm{U}_q^q>0.
\]
The eigenfunction equation and \(t_0^2=\lambda_1\) instead give the
identity \(-\Delta U=2t_0 M\partial_1\phi\).
Together these imply
\(2t_0\partial_1\phi=c M^{q-2}\phi^{q-1}>0\) almost everywhere,
contradicting \(\int_\Omega\partial_1\phi=0\), which follows by
density of \(C_c^\infty(\Omega)\) in \(H\). Thus \(g\ne0\).

Extend \(h_s\) by \eqref{eq:parameterization} to the closed unit ball
in \(X\). Define \(G_s(z)=R(e^{-tx_1}(\phi+h_s(z)))\). We use
derivatives in the \(B\) norm. The estimate
\(\norm{z}_2^2\le(\lambda_2-\lambda_1)^{-1}\norm{z}_B^2\)
makes the denominator uniformly positive for small \(s\), and
differentiation gives, uniformly on the ball,
\[
 h_s(z)=sz+O(s^3),\qquad Dh_s(z)=s I+O(s^3),
 \qquad D^2h_s(z)=O(s^3).
\]
Multiplication by \(e^{-tx_1}\) differs from multiplication by \(M\)
by \(O(s^2)\) in operator norm on \(H\), since \(t-t_0=O(s^2)\).
The map \(u\mapsto\int_\Omega |u|^q\) is twice continuously
Fr\'echet differentiable on \(H\) because \(q>2\) and
\(H\hookrightarrow L^q\) continuously. Its second derivative is
the bilinear form \(q(q-1)\int_\Omega|u|^{q-2}hk\), bounded by
H\"older's inequality. Thus \(D^2R\) is bounded on a sufficiently
small \(H\)-neighborhood of \(U\). The chain rule and Taylor's
formula give
\begin{equation}\label{eq:uniform-estimates}
 \begin{aligned}
 G_s(z)&=R(U)+s\ell(z)+O(s^2),\\
 \nabla_B G_s(z)&=sg+O(s^2),\qquad
 \norm{D^2G_s(z)}=O(s^2),
 \end{aligned}
\end{equation}
uniformly on the closed unit ball. Here \(\nabla_B\) denotes the
Riesz representative of the derivative with respect to \(B\).

For small \(s>0\) the normalized gradient map
\(T_s(z)=\nabla_B G_s(z)/\norm{\nabla_B G_s(z)}_B\)
is defined on the ball and takes values in its unit sphere. The
denominator is at least \(s\norm{g}_B/2\) and the Lipschitz constant
of the unnormalized gradient is \(O(s^2)\). Normalization therefore
makes \(T_s\) a contraction with Lipschitz constant \(O(s)\).
By the Banach fixed-point theorem, \(T_s\) has a fixed point \(z_s\),
which satisfies \(z_s=g/\norm{g}_B+O(s)\) by
\eqref{eq:uniform-estimates}. Set
\(\gamma_s=\norm{\nabla_B G_s(z_s)}_B\). For any \(z\) on the unit
sphere, \(B(z_s,z-z_s)=-\norm{z-z_s}_B^2/2\). Taylor's formula
on the segment inside the ball gives
\[
 G_s(z)-G_s(z_s)
 \le -\tfrac12\bigl(\gamma_s-O(s^2)\bigr)\norm{z-z_s}_B^2 .
\]
Since \(\gamma_s\ge s\norm{g}_B/2\), this proves strict global
maximality at \(z_s\). The first line of \eqref{eq:uniform-estimates}
then gives \eqref{eq:expansion} with \(K_\Omega=\norm{g}_B>0\).
\end{proof}

In terms of an \(L^2\)-orthonormal Dirichlet eigenbasis \((\phi_j)\),
with \(\phi_1=\phi\), \eqref{eq:riesz} gives the convergent series
\[
 K_\Omega^2=\sum_{j\ge2}\frac{|\ell(\phi_j)|^2}{\lambda_j-\lambda_1}.
\]
The normalized extremal \(u_t\), specified by
\(\inner{e^{tx_1}u_t}{\phi}_{L^2}=1\), satisfies
\[
 e^{tx_1}u_t=\phi+\sqrt{t^2-\lambda_1}\,g/\norm{g}_B
                      +O_H(t^2-\lambda_1)
\]
by \eqref{eq:parameterization}.

\begin{proposition}\label{prop:cylinder}
For \(\Omega=(0,L)\times D_0\), where \(L>0\) and
\(D_0\subset\R^{n-1}\) is a bounded connected Lipschitz domain,
there are positive constants \(c_1,c_2\) such that for all
sufficiently large \(t\),
\[
 c_1 e^{tL}t^{1-1/n}
 \le C_{\BJ}(e^{2tx_1},2,\Omega)
 \le c_2 e^{tL}t^{1-1/n}.
\]
\end{proposition}

\begin{proof}
Here \(\min_{\overline\Omega}a_t=1\) and \(\tfrac12\Delta a_t=2t^2a_t\).
For \(Q_{a_t}(u)=0\), \eqref{eq:form} and H\"older's inequality give
\[
 E(u)\le 2t^2\int_\Omega a_t u^2
 \le 2t^2\left(|D_0|\frac{e^{ntL}-1}{nt}\right)^{2/n}\norm{u}_q^2.
\]
For the lower estimate, choose a nonzero
\(\psi\in C_c^\infty(1,K)\), with \(K>1\) large enough that
\(\norm{\psi'}_2^2<\norm{\psi}_2^2\). Such a cutoff is obtained
by dilation of a fixed smooth bump. Choose a nonzero
\(\eta\in C_c^\infty(D_0)\). Let
\(v(x)=\psi(t(L-x_1))\eta(x')\) and \(u=e^{-tx_1}v\).
The support is inside the cylinder for large \(t\), and
\eqref{eq:gauge} gives
\[
 Q_{a_t}(u)
 =t(\norm{\psi'}_2^2-\norm{\psi}_2^2)\norm{\eta}_2^2
   +t^{-1}\norm{\psi}_2^2\norm{\nabla\eta}_2^2\le-ct .
\]
Changing variables in the \(L^q\) norm gives
\(\norm{u}_q=e^{-tL}t^{-1/q}\norm{e^y\psi(y)}_q\norm{\eta}_q\).
In \eqref{eq:relaxation}, the numerator at \(u\) is at least
\(-Q_{a_t}(u)\ge ct\), giving the claim by \(1/2+1/q=1-1/n\).
\end{proof}

\section{Sharp reverse estimates in \texorpdfstring{\(L^p\)}{Lp}}\label{sec:general-lp}

Let \(n\ge2\), \(1<p<n\), and \(p^*=np/(n-p)\). We write \(D\)
for either a bounded Lipschitz domain or \(\R^n\) and set
\(X_p(D)=W^{1,p}_0(D)\) in the former case and
\(X_p(\R^n)=\dot W^{1,p}(\R^n)\) in the latter. Suppose \(a\) is
bounded and continuous on \(D\), with \(m=\inf_D a>0\) and weak
gradient \(\nabla a\in L^n(D)\). The product rule and the Sobolev
embedding give \(au\in X_p(D)\) for \(u\in X_p(D)\).
Define
\begin{equation}\label{eq:general-lp-pairing}
 Q_{a,p}(u)=\int_D|\nabla u|^{p-2}\nabla u
                         \cdot\nabla(au)\,dx .
\end{equation}
Here \(\abs{\xi}^{p-2}\xi\) is set to zero at \(\xi=0\).
By convexity and differentiability of the \(p\)-th power of the
\(L^p\) norm, \(Q_{a,p}(u)=0\) is equivalent to
\(\nabla u\perp_{\BJ}\nabla(au)\) in \(L^p(D;\R^n)\). If there is
a nonzero orthogonal profile, put
\[
 C_{\BJ,p}(a,D)=
 \sup_{\substack{0\ne u\in X_p(D)\\Q_{a,p}(u)=0}}
                  \frac{\norm{\nabla u}_{L^p(D)}}{\norm{u}_{L^{p^*}(D)}} .
\]

\begin{theorem}\label{thm:general-lp-relaxation}
For the weights and spaces above, if \(C_{\BJ,p}(a,D)\) is defined,
then
\begin{equation}\label{eq:general-lp-relaxation}
 \displaystyle
 C_{\BJ,p}(a,D)^p
 =
 \sup_{\substack{0\ne u\in X_p(D)\\ Q_{a,p}(u)\le0}}
 \frac{\displaystyle \int_D|\nabla u|^p\,dx
       -m^{-1}\int_D|\nabla u|^{p-2}\nabla u\cdot\nabla(au)\,dx}
      {\displaystyle\left(\int_D|u|^{np/(n-p)}\,dx\right)^{(n-p)/n}} .
\end{equation}
The orthogonality supremum is attained if and only if the relaxed
problem on the right has a maximizer \(u\) satisfying \(Q_{a,p}(u)=0\).
If a nonzero \(u\in X_p(D)\) satisfies \(Q_{a,p}(u)\le0\),
then the equality constraint admits a nonzero profile.
\end{theorem}

\begin{proof}
The product rule gives
\[
 Q_{a,p}(u)=\int_D a|\nabla u|^p\,dx
 +\int_D u|\nabla u|^{p-2}\nabla u\cdot\nabla a\,dx .
\]
For \(Q_{a,p}(u)\le0\), H\"older's inequality yields
\(m\norm{\nabla u}_p\le \norm{\nabla a}_n\norm u_{p^*}\).
In particular, the orthogonality supremum is finite when defined.

Fix \(0\ne u\in X_p(D)\) with \(Q_{a,p}(u)<0\). Choose balls
\(B_j\Subset D\) with radii tending to zero and
\(\sup_{B_j}(a-m)\to0\). Continuity and the definition of \(m\)
give such balls, even when their centers tend to infinity.
Take \(w_j\in C_c^\infty(B_j)\) with
\(\norm{\nabla w_j}_p=1\) and \(\norm{w_j}_{p^*}\le 1/j\).
For each ball one can multiply a nonzero smooth cutoff by
\(\sin(k x_1)\) and normalize its gradient. The leading part of
the gradient has \(L^p\) norm of order \(k\) as \(k\to\infty\),
by periodic averaging, whereas the unnormalized \(L^{p^*}\)
norm remains bounded.

Uniformly for \(s\) in each bounded subset of \([0,\infty)\),
\begin{align}
 Q_{a,p}(u+s w_j)&=Q_{a,p}(u)+ms^p+o(1),\label{eq:general-lp-recover-q}\\
 \int_D|\nabla(u+s w_j)|^p\,dx
   &=\int_D|\nabla u|^p\,dx+s^p+o(1).\label{eq:general-lp-recover-e}
\end{align}
To see the gradient estimates, the integrals outside \(B_j\)
are unchanged, and \(\norm{\nabla u}_{L^p(B_j)}\to0\).
On \(B_j\), the \(L^p\) norm of \(\nabla(u+s w_j)\) therefore
differs from \(s\) by \(o(1)\) by the reverse triangle inequality.
The analogous weighted norm differs from
\(s\norm{a^{1/p}\nabla w_j}_p\) by \(o(1)\), and
\(\norm{a^{1/p}\nabla w_j}_p^p\to m\).
The mixed integral on \(B_j\) is bounded in absolute value by
\[
 \norm{\nabla a}_{L^n(B_j)}\,\norm{u+s w_j}_{L^{p^*}(B_j)}
                      \norm{\nabla(u+s w_j)}_{L^p(B_j)}^{p-1},
\]
which tends to zero uniformly for bounded \(s\). The same estimate
with \(s=0\) disposes of the original mixed integral on \(B_j\).
This proves \eqref{eq:general-lp-recover-q} and
\eqref{eq:general-lp-recover-e}.

Continuity in \(s\) and \eqref{eq:general-lp-recover-q} give roots
\(s_j\ge0\) in a fixed bounded interval, with
\[
 Q_{a,p}(u+s_jw_j)=0,\qquad
 s_j\longrightarrow\bigl(-Q_{a,p}(u)/m\bigr)^{1/p}.
\]
The \(L^{p^*}\) norms tend to \(\norm u_{p^*}>0\), and
\eqref{eq:general-lp-recover-e} gives the numerator on the right of
\eqref{eq:general-lp-relaxation} in the limit. For \(Q_{a,p}(u)=0\)
that numerator equals \(\norm{\nabla u}_p^p\) directly.
Taking the respective suprema proves the formula and the
feasibility assertion. On \(Q_{a,p}(u)=0\), the relaxed quotient
equals the \(p\)-th power of the original quotient. Hence each
extremal for the original problem attains the relaxed supremum,
and each relaxed maximizer satisfying \(Q_{a,p}(u)=0\) attains
the original supremum.
\end{proof}

\begin{proposition}\label{prop:general-lp-cylinder}
Let \(D=(0,L)\times D_0\), with \(L>0\) and \(D_0\subset\R^{n-1}\)
a bounded connected Lipschitz domain. For every \(1<p<n\),
there exist constants \(c,C>0\) such that for all sufficiently
large \(t\), the weight \(a_t(x)=e^{ptx_1}\) admits a nonzero
orthogonal profile and
\begin{equation}\label{eq:general-lp-cylinder}
 \displaystyle
   c\,e^{tL}t^{1-1/n}
   \le C_{\BJ,p}(e^{ptx_1},D)
   \le C\,e^{tL}t^{1-1/n}.
\end{equation}
\end{proposition}

\begin{proof}
For \(Q_{a_t,p}(u)=0\), the identity \(\nabla a_t=pt a_t e_1\)
and weighted H\"older give
\[
 \int_D a_t|\nabla u|^p\,dx
 \le pt\left(\int_D a_t|\nabla u|^p\,dx\right)^{(p-1)/p}
        \left(\int_D a_t|u|^p\,dx\right)^{1/p}.
\]
Since \(\inf_D a_t=1\), another application of H\"older gives
\[
 \norm{\nabla u}_p^p
 \le (pt)^p\int_D a_t|u|^p\,dx
 \le (pt)^p
     \left(\frac{|D_0|(e^{ntL}-1)}{nt}\right)^{p/n}
     \norm u_{p^*}^p.
\]
This gives the upper bound.

For the lower bound, choose nonzero \(\eta\in C_c^\infty(D_0)\).
For \(\psi\in C_c^\infty((0,\infty))\) to be fixed below, put
\[
 y=t(L-x_1),\qquad v_t(x)=\psi(y)\eta(x'),\qquad
 u_t(x)=e^{-tx_1}v_t(x).
\]
The support lies in \(D\) for sufficiently large \(t\). The
exponential factors cancel in the pairing, giving
\[
 Q_{a_t,p}(u_t)
 =\int_D|\nabla v_t-t v_te_1|^{p-2}(\nabla v_t-t v_te_1)
                           \cdot(\nabla v_t+(p-1)t v_te_1)\,dx .
\]
Changing variables to \(y\) and dividing the gradient factors
by \(t\), dominated convergence gives
\begin{equation}\label{eq:general-lp-cylinder-limit}
 \lim_{t\to\infty}t^{1-p}Q_{a_t,p}(u_t)
  =\norm{\eta}_{L^p(D_0)}^p
   \int_0^\infty|\psi'+\psi|^{p-2}(\psi'+\psi)
                         \bigl(\psi'-(p-1)\psi\bigr)\,dy .
\end{equation}
The transverse derivatives divided by \(t\) tend to zero on a
fixed compact support in \((y,x')\), so the argument also holds
for \(1<p<2\) by continuity of \(\xi\mapsto|\xi|^{p-2}\xi\).
Take \(\psi(y)=\chi(y/R)\), where \(0\ne\chi\in C_c^\infty(1,2)\).
After division by \(R\) and rescaling \(y\), the integral on the
right of \eqref{eq:general-lp-cylinder-limit} tends to
\(-(p-1)\int_1^2|\chi|^p\,dy\) as \(R\to\infty\). Fix \(R\)
large enough that this integral before division is negative.
Then \(Q_{a_t,p}(u_t)\le-c_0t^{p-1}\) for large \(t\), so
Theorem~\ref{thm:general-lp-relaxation} gives feasibility. Also,
\[
 \norm{u_t}_{p^*}
 =e^{-tL}t^{-1/p^*}\norm{e^y\psi(y)}_{L^{p^*}(0,\infty)}
                         \norm{\eta}_{L^{p^*}(D_0)} .
\]
Using \(u_t\) in \eqref{eq:general-lp-relaxation} with \(m=1\)
gives the lower bound since \((p-1)/p+1/p^*=1-1/n\).
\end{proof}

\begin{proposition}\label{prop:general-lp-weak-lsc}
If \(1<p<n\) and \(p\ne2\), every nonempty bounded Lipschitz
domain \(D\) admits a smooth positive bounded weight \(a\) for
which \(Q_{a,p}\) fails to be sequentially weakly lower
semicontinuous on \(W^{1,p}_0(D)\).
\end{proposition}

\begin{proof}
Take \(a(x)=A+x_1>0\) on \(\overline D\). Fix a ball \(B\Subset D\)
and a smooth compactly supported \(u\) which is affine on \(B\),
with \(u\ge1/2\) there and \(\nabla u=\sigma e_1\), where
\(\sigma=-\operatorname{sgn}(p-2)\delta\) for a small \(\delta>0\).
For a gradient \(s e_1\), the integrand at \(u(x)\) is
\(g_x(s)=a(x)|s|^p+u(x)|s|^{p-2}s\). Its second derivative
for \(s\ne0\) is
\[
 g_x''(s)=(p-1)|s|^{p-3}
          \bigl(p a(x)|s|+(p-2)u(x)\operatorname{sgn}(s)\bigr).
\]
For sufficiently small \(\delta\), there is a uniform
\(\kappa>0\) with \(g_x''(s)\le-\kappa\) on \(B\) whenever
\(|s-\sigma|\le\delta/2\).
Choose \(0\le\chi\le1\) smooth, nonzero and compactly supported
in \(B\), and set
\[
 u_j(x)=u(x)+\frac{\varepsilon}{j}\chi(x)\sin(jx_1),
 \qquad 0<\varepsilon<\delta/2.
\]
Then \(u_j\to u\) uniformly, and \(u_j\rightharpoonup u\) in
\(W^{1,p}_0(D)\) by periodic averaging. On \(B\) the gradient
equals \((\sigma+\varepsilon\chi\cos(jx_1))e_1+O(j^{-1})\)
uniformly. Continuity on bounded sets lets us replace the
integrand there by \(g_x(\sigma+\varepsilon\chi\cos(jx_1))\)
with an \(o(1)\) error in the integral. Taylor's inequality now
gives
\[
 Q_{a,p}(u_j)-Q_{a,p}(u)
 \le \varepsilon\int_B g_x'(\sigma)\chi\cos(jx_1)\,dx
       -\frac{\kappa\varepsilon^2}{2}\int_B\chi^2\cos^2(jx_1)\,dx
       +o(1).
\]
The Riemann--Lebesgue lemma sends the first integral to zero and
the second to \(\tfrac12\int_B\chi^2\,dx>0\). The limit superior
is strictly negative.
\end{proof}

Existence of a relaxed maximizer in the quadratic case, as proved in
Theorems~\ref{thm:relaxation} and~\ref{thm:whole-space-relaxation},
uses weak lower semicontinuity of the pairing.

\section{Examples and applications}\label{sec:pde}

\noindent\textbf{The sharp constant on a four-dimensional cube.}
For \(\Omega=(0,\pi)^4\), we have \(t_0=2\), \(q=4\), and
\(U=e^{-2x_1}\prod_{j=1}^4\sin x_j\) up to a constant. Separating
the integrals in \eqref{eq:threshold-value} and using
\[
 \int_0^\pi e^{-4x}\sin^2 x\,dx=\frac{1-e^{-4\pi}}{40},
 \qquad
 \int_0^\pi e^{-8x}\sin^4 x\,dx=\frac{3(1-e^{-8\pi})}{5440},
\]
together with the unweighted integrals \(\pi/2\) and \(3\pi/8\),
respectively, yields
\begin{equation}\label{eq:cube}
 C_{\BJ}\bigl(e^{4x_1},2,(0,\pi)^4\bigr)
 =\left(\frac{8704}{405}\pi^3\tanh(2\pi)\right)^{1/4}.
\end{equation}

\smallskip
\smallskip
\noindent\textbf{A radial threshold on the whole space.}
Let \(n\ge3\), \(r=|x|\), and \(a_t(x)=\exp(2t\arctan r^2)\) for \(t\ge0\).
This is a smooth weight with \(1\le a_t\le e^{\pi t}\) and
\(\nabla a_t\in L^n\). Put \(g=\arctan r^2\) and \(v=e^{tg}u\).
Multiplication by \(e^{\pm tg}\) is bounded on \(\dot H^1\), and
expansion of the gradients gives
\begin{equation}\label{eq:whole-radial-transform}
 \int_{\R^n}\nabla u\cdot\nabla(a_tu)\,dx
 =\int_{\R^n}|\nabla v|^2\,dx
   -4t^2\int_{\R^n}\frac{|x|^2}{(1+|x|^4)^2}\,v^2\,dx.
\end{equation}
The threshold is given by
\begin{equation}\label{eq:whole-radial-threshold}
 \displaystyle
 \exists\,0\ne u\in\dot H^1(\R^n),\quad
 \nabla u\perp_{\BJ}\nabla(a_tu)
 \quad\Longleftrightarrow\quad t\ge t_0=\frac{\sqrt{n^2-4}}{2}.
\end{equation}
The normalized positive radial solution of the critical H\'enon equation
with weight exponent \(2\) is \(v_0=(1+r^4)^{-(n-2)/4}\)
\cite[Eq.~(1.8)]{GladialiGrossiNeves2013}. It belongs to \(\dot H^1\)
and satisfies
\[
 -\Delta v_0=(n^2-4)\frac{|x|^2}{(1+|x|^4)^2}\,v_0 .
\]
For compactly supported smooth \(v\), integration by parts gives
\[
 \int_{\R^n}|\nabla v|^2\,dx
 -(n^2-4)\int_{\R^n}\frac{|x|^2}{(1+|x|^4)^2}\,v^2\,dx
 =\int_{\R^n}v_0^2|\nabla(v/v_0)|^2\,dx.
\]
The left side is continuous on \(\dot H^1\) by H\"older and Sobolev,
since \(r^2(1+r^4)^{-2}\in L^{n/2}\). It is therefore nonnegative
on \(\dot H^1\). If it vanishes, approximation and local \(H^1\)
convergence imply \(\nabla(v/v_0)=0\) on each bounded ball. Testing
the equation for \(v_0\) with \(v_0\) by density shows that \(v_0\)
gives equality. Thus at \(t=t_0\) the nonzero orthogonal profiles
are exactly the nonzero multiples of
\[
 u_0(x)=e^{-t_0\arctan |x|^2}(1+|x|^4)^{-(n-2)/4}.
\]
For \(t<t_0\) the form in \eqref{eq:whole-radial-transform} is
strictly positive on nonzero \(v\). For \(t>t_0\), \(v_0\) gives a
negative direction, while rapid oscillation of a compactly supported
cutoff gives a positive direction. Continuity on their span gives
a nonzero vector on which the form vanishes. This proves
\eqref{eq:whole-radial-threshold}. At the threshold,
\(C_{\BJ}(a_{t_0},2,\R^n)=\norm{\nabla u_0}_2/\norm{u_0}_{2n/(n-2)}\).
Since \(u_0\) decays as \(r^{2-n}\), it is in \(H^1\) exactly when
\(n\ge5\). In \(n=3,4\) there is no nonzero \(H^1\) profile at
the threshold. Above the threshold, approximation of \(v_0\) gives a
compactly supported negative direction, so there is also a nonzero
orthogonal profile in \(H^1\).

\smallskip
\noindent\textbf{A periodic weight with infinite best constant.}
On \(\R^n\), \(n\ge3\), take \(a(x)=\exp(2\sin x_1)\).
For nonzero \(\chi\in C_c^\infty(\R^n)\) let
\(u_R=a^{-1/2}\chi(x/R)\). Expansion as in
\eqref{eq:whole-radial-transform} gives
\[
 \int_{\R^n}\nabla u_R\cdot\nabla(au_R)\,dx
 =R^{n-2}\int_{\R^n}|\nabla\chi(y)|^2\,dy
  -R^n\int_{\R^n}\cos^2(Ry_1)\chi(y)^2\,dy .
\]
This expression is positive for small \(R\) and negative for
large \(R\), since the last integral tends to
\(\tfrac12\int\chi^2\) by the Riemann--Lebesgue lemma.
Continuity gives a nonzero compactly supported orthogonal profile
\(u_*\). Translating it along the \(x_2\)-direction preserves both
the pairing and the norms. The sum of \(N\) translates with disjoint
supports still satisfies the orthogonality constraint and has quotient
\[
 N^{1/n}\frac{\norm{\nabla u_*}_{L^2(\R^n)}}
                {\norm{u_*}_{L^{2n/(n-2)}(\R^n)}} .
\]
Hence \(C_{\BJ}(e^{2\sin x_1},2,\R^n)=\infty\).

\smallskip
\noindent\textbf{The Dirichlet heat equation.}
Testing \(\partial_\tau u=\Delta u\) with \(au\) gives
\begin{equation}\label{eq:heat}
 \frac12\frac{d}{d\tau}\int_\Omega a u(\tau)^2=-Q_a(u(\tau)).
\end{equation}
The use of exponential weights in heat semigroup estimates is
classical \cite{MuruganSaloffCoste2017}. For \(a=e^{2tx_1}\),
the identity \eqref{eq:heat} and the substitution in \eqref{eq:gauge} give
\[
 \frac{\dfrac{d}{d\tau}\int_\Omega a u^2}{2\int_\Omega a u^2}
 =t^2-\frac{E(e^{tx_1}u)}{\norm{e^{tx_1}u}_2^2}
 \le t^2-\lambda_1 .
\]
The initial datum \(u(0)=e^{-tx_1}\phi\) attains the upper bound.
If \(t^2>\lambda_1\), continuity of the heat flow in \(H\) gives strict
initial growth. The Dirichlet spectral gap and boundedness of \(a\)
give decay to zero as \(\tau\to\infty\). Thus the weighted square norm
has a positive interior maximum, where the nonzero profile satisfies
\(Q_a(u)=0\).

\smallskip
\noindent\textbf{The nonlinear \(p\)-heat equation.}
For \(\partial_\tau u=\Delta_p u\),
testing with \(au\) in the energy identity gives
\begin{equation}\label{eq:pheat}
 -\frac12\frac{d}{d\tau}\int_\Omega a u^2
 =\int_\Omega|\nabla u|^{p-2}\nabla u\cdot\nabla(au).
\end{equation}
For \(1<p<\infty\), the vanishing of the right-hand side is
equivalent to the Birkhoff--James condition in \(L^p\), by
convexity and differentiability of the \(p\)-th power of the norm.
For an explicit example take \(2<p<n\). The compactly supported
Barenblatt solution has the form
\(u(\tau,x)=\tau^{-n\beta}F(\tau^{-\beta}x)\), where
\(\beta=(n(p-2)+p)^{-1}\) and
\(F(y)=(c-\kappa |y|^{p/(p-1)})_+^{(p-1)/(p-2)}\) for \(c>0\) and
\(\kappa=\frac{p-2}{p}\beta^{1/(p-1)}\)
\cite{KaminVazquez1988}; see also \cite[Chap.~11]{Vazquez2006}.
For an even integer \(k>n\), use \(a(x)=1+|x|^k\) and set
\(A=\int_{\R^n}F^2\), \(B_0=\int_{\R^n}|y|^kF(y)^2\,dy\).
The weighted square norm is
\(A\tau^{-n\beta}+B_0\tau^{(k-n)\beta}\), whose derivative
vanishes exactly at
\[
 \tau_*=\left(\frac{nA}{(k-n)B_0}\right)^{1/(k\beta)} .
\]
One may work on a fixed ball containing the supports for times in a
neighborhood of \(\tau_*\). There the spatial trace is zero and the
weight is smooth, bounded, and positive. The explicit profile is
\(C^1\), including at the origin and at its free boundary. Both
\(\partial_\tau u\) and \(\nabla u\) are continuous on the ball for
positive times. The weak identity \eqref{eq:pheat}, first obtained
in distributions in time by testing with \(au\), therefore holds
pointwise in time in this neighborhood.

\smallskip
\noindent\textbf{A two-component convection--diffusion system.}
At the exponential threshold \(t_0=\sqrt{\lambda_1}\), the extremal
\(U=e^{-t_0x_1}\phi\) satisfies
\[
 -\Delta U-2t_0\partial_1 U=2t_0^2 U .
\]
This also yields an explicit two-component evolution on the same
domain. Take \(J=\left(\begin{smallmatrix}0&-1\\1&0\end{smallmatrix}\right)\)
and \(\omega\in\R\). The function
\[
 W(\tau,x)= e^{-2t_0^2\tau}U(x)
        \bigl(\cos(\omega\tau),\sin(\omega\tau)\bigr)^{\mathsf T}
\]
solves \(\partial_\tau W=\Delta W+2t_0\partial_1W+\omega JW\)
with zero Dirichlet boundary values. Its components satisfy
\(Q_{a_{t_0}}(W_j)=0\) for every \(\tau\). Each nonzero component
attains the scalar constant \eqref{eq:threshold-value}.

\smallskip
\noindent\textbf{The incompressible Navier--Stokes system.}
On the whole space or on a periodic domain, the fluid applications
consider the unforced system for the velocity \(u\) and pressure \(p\),
with kinematic viscosity \(\nu>0\) and divergence-free initial data \(u_0\),
\begin{equation}\label{eq:ns-system}
 \begin{cases}
 \partial_t u+(u\cdot\nabla)u+\nabla p=\nu\Delta u,\\
 \nabla\cdot u=0,\\
 u(0,x)=u_0(x).
 \end{cases}
\end{equation}
The cylindrical example below has an axially periodic domain and a
no-slip lateral wall. For smooth finite-energy solutions on \(\R^3\),
testing \eqref{eq:ns-system} with \(u\) gives the kinetic energy identity
\begin{equation}\label{eq:ns-kinetic}
 \frac12\|u(t)\|_{L^2(\R^3)}^2
  +\nu\int_0^t\|\nabla u(s)\|_{L^2(\R^3)}^2\,ds
  =\frac12\|u_0\|_{L^2(\R^3)}^2 .
\end{equation}
The corresponding energy inequality underlies Leray's weak solution
theory \cite{Leray1934}. Classical interior regularity criteria use
space--time integrability of \(u\) \cite{Serrin1962}. The scale-critical
endpoint \(L^\infty_tL^3_x\) in three dimensions is treated in
\cite{EscauriazaSereginSverak2003}. We use reverse estimates to
control the stronger spatial norms needed for continuation.

\smallskip
\noindent\textbf{Axisymmetric Navier--Stokes flow with swirl.}
We work on \(\R^3\) and use cylindrical components
\(u=u^r e_r+u^\theta e_\theta+u^z e_z\), with
\(r=(x_1^2+x_2^2)^{1/2}\) and \(z=x_3\). Let
\(F=u^\theta/r\) and \(\Gamma=ru^\theta\).
For \(r>0\), the azimuthal component of \eqref{eq:ns-system} and
the equation for \(\Gamma\) take the form
\begin{equation}\label{eq:ns-swirl-system}
 \begin{aligned}
 \partial_t u^\theta+u^r\partial_r u^\theta
   +u^z\partial_z u^\theta+\frac{u^r}{r}u^\theta
 &=\nu\left(\partial_r^2+\frac1r\partial_r
             +\partial_z^2-\frac1{r^2}\right)u^\theta,\\
 \partial_t\Gamma+u^r\partial_r\Gamma+u^z\partial_z\Gamma
 &=\nu\left(\partial_r^2-\frac1r\partial_r+\partial_z^2\right)\Gamma .
 \end{aligned}
\end{equation}
The second equation gives the maximum principle
\(\norm{\Gamma(t)}_\infty\le\norm{\Gamma(0)}_\infty\)
\cite[Eq.~(2.2)]{Palasek2026}. The analysis of possible axisymmetric
singularities uses Liouville theorems for rescaled flows
\cite{KochNadirashviliSereginSverak2009}, and azimuthal velocity
integrability gives one-component regularity criteria
\cite{ChenFangZhang2017}. Our continuation argument uses the
circulation bound and the gradient pairing to verify such an
integrability criterion.

\begin{theorem}\label{thm:axisym-continuation}
Let \(u\) be a smooth axisymmetric solution of the unforced
Navier--Stokes equations on \([0,T_*)\times\R^3\), with viscosity
\(\nu>0\), \(T_*<\infty\), \(u\in C([0,T_*);H^3(\R^3))\), and
\(ru^\theta(0)\in L^\infty(\R^3)\). Fix
\(a\in L^\infty(\R^3)\) with
\(\operatorname*{ess\,inf}_{\R^3}a>0\) and weak gradient
\(\nabla a\in L^3(\R^3)\). If, for some \(0\le t_0<T_*\),
\begin{equation}\label{eq:axisym-defect}
 \int_{t_0}^{T_*}
 \left\|\frac{u^\theta(t)}r\right\|_{L^2(\R^3)}
 \left[
   \int_{\R^3}\nabla\!\left(\frac{u^\theta(t)}r\right)
     \cdot\nabla\!\left(a\,\frac{u^\theta(t)}r\right)\,dx
 \right]_+^{1/2}\,dt<\infty,
\end{equation}
where \([s]_+=\max(s,0)\), then \(u\) extends smoothly beyond \(T_*\).
In particular, it suffices that
\(\nabla(u^\theta/r)\perp_{\BJ}\nabla(a u^\theta/r)\) in
\(L^2(\R^3;\R^3)\) for almost every \(t\in(t_0,T_*)\).
\end{theorem}

\begin{proof}
At \((r,0,z)\), axisymmetry gives \(F=-\partial_{x_2}u_1\).
Differentiating this identity in \(r\) and \(z\), then rotating about
the axis, bounds \(|F|\) and \(|\nabla F|\) by \(|\nabla u|\) and
\(|D^2u|\), respectively. In particular \(F\in
C([0,T_*);H^1(\R^3))\). The axis has measure zero, and the smooth
axisymmetric quotient extends across it. Since \(a\) is bounded and
\(\nabla a\in L^3\), the Sobolev embedding gives \(aF\in H^1\).

All norms in this proof are on \(\R^3\). Put
\(m=\operatorname*{ess\,inf}a\) and
\(Q_a(F)=\int_{\R^3}\nabla F\cdot\nabla(aF)\,dx\). By the product rule,
\[
 m\norm{\nabla F}_2^2
 \le [Q_a(F)]_+ +\norm{\nabla F}_2\norm{F\nabla a}_2 .
\]
Let \(S\) be a Sobolev constant for \(\norm{F}_6\le S\norm{\nabla F}_2\).
Choose \(K\) with \(S\norm{\nabla a\,
 \mathbf 1_{\{|\nabla a|>K\}}}_3\le m/2\).
H\"older's inequality on the two parts of \(\nabla a\) gives
\(\norm{F\nabla a}_2\le K\norm F_2+(m/2)\norm{\nabla F}_2\).
Solving the resulting quadratic inequality yields
\begin{equation}\label{eq:axisym-multiplier}
 \norm{\nabla F}_2\le \frac{2K}{m}\norm F_2+
                 \sqrt{\frac{2[Q_a(F)]_+}{m}} .
\end{equation}
The angular derivative of the velocity satisfies
\(r^{-1}\partial_\theta u=(u^r/r)e_\theta-(u^\theta/r)e_r\).
Combining the angular derivative with \eqref{eq:ns-kinetic} gives
\[
 \int_0^{T_*}\norm{F(t)}_2^2\,dt
 \le\int_0^{T_*}\norm{\nabla u(t)}_2^2\,dt
 \le \frac{\norm{u(0)}_2^2}{2\nu}.
\]
Now \((u^\theta)^2=\Gamma F\). Using the circulation maximum
principle, interpolation, and the Sobolev inequality, we obtain
\[
 \norm{u^\theta(t)}_6^4
 \le \norm{\Gamma(0)}_\infty^2\norm{F(t)}_3^2
 \le S\norm{\Gamma(0)}_\infty^2\norm{F(t)}_2\norm{\nabla F(t)}_2 .
\]
Integrating and using \eqref{eq:axisym-defect}
and \eqref{eq:axisym-multiplier} gives
\(u^\theta\in L^4((t_0,T_*);L^6(\R^3))\).
On earlier compact time intervals the same norm is finite by the
assumed Sobolev regularity. The single-component continuation
criterion of Zhang and Zhang for \(u^\theta\in L^4_tL^6_x\) gives
the smooth extension \cite{ZhangZhang2014}. The same estimates give
\(u^\theta/r\in L^2_tL^3_x\) by interpolation, which is the \(d=1\)
case of the weighted criterion of Chen, Fang, and Zhang
\cite[Remark~1.3]{ChenFangZhang2017}.
\end{proof}

The same proof gives a criterion for weights chosen separately at each
time. Under the flow hypotheses of Theorem~\ref{thm:axisym-continuation},
it suffices to have a measurable \(K(t)\ge0\) on \((t_0,T_*)\) with
\begin{equation}\label{eq:swirl-lower-order}
 \left\|\nabla\!\left(\frac{u^\theta(t)}r\right)\right\|_2
 \le K(t)\left\|\frac{u^\theta(t)}r\right\|_2,\qquad
 \int_{t_0}^{T_*}K(t)\left\|\frac{u^\theta(t)}r\right\|_2^2\,dt<\infty .
\end{equation}
Indeed, the last estimate for \(\norm{u^\theta(t)}_6^4\) in the proof
is integrable under \eqref{eq:swirl-lower-order}. The first inequality
holds with \(K(t)=\norm{\nabla a_t}_\infty/\inf a_t\) whenever
\(a_t\in W^{1,\infty}\), \(\inf a_t>0\) and
\(\int\nabla F(t)\cdot\nabla(a_t F(t))\le0\), by the product rule.
A uniform bound on \(K(t)\) then suffices by the kinetic energy estimate.

For example, if the initial circulation \(ru^\theta(0)\) has one sign,
the circulation equation preserves that sign, and so does
\(F=\Gamma/r^2\) away from the axis. Proposition~\ref{prop:profile-weight}
therefore supplies an axisymmetric positive weight at every time with
a nonzero \(C^2\cap H^2\) profile. A quantitative spatial condition
ensuring continuation is the existence of a single nonnegative
\(\chi\in C_c^\infty(\R^3)\) and \(\delta>0\) such that
\[
 \frac12\int_{\R^3}(\Delta\chi)F(t)^2\,dx
       -\int_{\R^3}\chi|\nabla F(t)|^2\,dx
 \ \ge\ \delta\int_{\R^3}|\nabla F(t)|^2\,dx
\]
for almost every \(t\in(t_0,T_*)\). In fact \(a=1+\chi/\delta\)
then satisfies \(\int\nabla F(t)\cdot\nabla(aF(t))\le0\) by
integration by parts. The same cutoff and the same \(\delta\) control
the full terminal interval.

\begin{proposition}\label{prop:axisym-obstruction}
There are smooth, divergence-free, axisymmetric initial data of finite
energy and bounded circulation for which the weighted gradient pairing is
strictly positive for every weight in Theorem~\ref{thm:axisym-continuation}.
One such datum is
\[
 F_0(x)=\frac{2z^2-r^2}{(1+|x|^2)^{5/2}},\qquad
 u_0(x)=(-x_2,x_1,0)F_0(x).
\]
For this datum,
\begin{equation}\label{eq:axisym-positive-pairing}
 \int_{\R^3}\nabla F_0\cdot\nabla(aF_0)\,dx
 =35\int_{\R^3}\frac{a(x)F_0(x)^2}{(1+|x|^2)^2}\,dx
 \ge(\operatorname*{ess\,inf}a)\int_{\R^3}|\nabla F_0|^2\,dx>0 .
\end{equation}
\end{proposition}

\begin{proof}
The field is pure swirl with \(u_0^\theta=rF_0\), and its divergence
vanishes since \(F_0\) is axisymmetric. The decay
\(F_0=O(|x|^{-3})\) gives \(u_0\in L^2\) and
\(ru_0^\theta=r^2F_0\in L^\infty\). The formula is smooth and gives
\(u_0\in H^k\) for every nonnegative integer \(k\).
Up to a constant factor, \(F_0\) is the conformal image under
stereographic projection of the degree-two spherical harmonic
\(\xi\mapsto 2\xi_3^2-\xi_1^2-\xi_2^2\) on \(\mathbb S^3\)
\cite[Eqs.~(2.1)--(2.2)]{Konig2023}.
With \(P=2z^2-r^2\) and \(s=|x|^2\), homogeneity and harmonicity give
\(\Delta(Pg(s))=P(14g'(s)+4s g''(s))\).
For \(g(s)=(1+s)^{-5/2}\), this yields
\(-\Delta F_0=35(1+|x|^2)^{-2}F_0\). Integrating by parts against
\(aF_0\in H^1\) proves the equality in
\eqref{eq:axisym-positive-pairing}. The lower bound follows by
testing the same Laplacian identity against \(F_0\).
\end{proof}

For each fixed admissible \(a\), the pairing is continuous on \(H^1\).
It therefore stays positive on a short interval along the local
smooth solution from Proposition~\ref{prop:axisym-obstruction}.
If a solution in the class of Theorem~\ref{thm:axisym-continuation}
has a finite first singular time, the integral in
\eqref{eq:axisym-defect} is infinite on every terminal interval
for each fixed admissible \(a\).

Substituting \(\Gamma=r^2 F\) into the circulation
equation in \eqref{eq:ns-swirl-system} gives
\[
 (\partial_t+u^r\partial_r+u^z\partial_z)F
 =\nu\left(\Delta+\frac2r\partial_r\right)F
        -2\frac{u^r}{r}F .
\]
For a smooth time-independent axisymmetric \(a=a(r,z)\),
testing with \(aF\) and integrating the transport term by parts
combines the transport and stretching contributions into
\[
 \frac12\int_{\R^3}
 (u^r\partial_r a+u^z\partial_z a-4a u^r/r)F^2\,dx .
\]
Universal pointwise cancellation forces
\(\partial_z a=0\) and \(\partial_r a=4a/r\), since a smooth local
stream function can prescribe \(u^r\) and \(u^z\) independently
at a point with \(r>0\). Hence \(a=cr^4\).
Its weighted mass is \(c\int_{\R^3}\Gamma^2\,dx\), and the loss of
the positive lower bound on the axis prevents the absorption used
in \eqref{eq:axisym-multiplier}.

There is also an exact flow for which a single positive weight works
at all times. On the axially periodic cylinder
\(D=B_1^2\times(\R/2\pi\mathbb Z)\) with a no-slip lateral wall, let
\(k\) be the first positive zero of the Bessel function \(J_1\). Set
\[
 u^r=u^z=0,\qquad u^\theta(t,r)=e^{-\nu k^2t}J_1(kr),\qquad
 \partial_r p=(u^\theta)^2/r,\quad \partial_zp=0 .
\]
The Bessel equation gives
\((\partial_r^2+r^{-1}\partial_r-r^{-2})J_1(kr)=-k^2J_1(kr)\);
the pressure balances the radial acceleration. Here
\(F(t,r)=e^{-\nu k^2t}\phi(r)\), where \(\phi(r)=J_1(kr)/r\)
extends smoothly across the axis. The scalar Laplacian satisfies
\(-\Delta\phi=k^2\phi+2\phi'/r\). Since \(\phi>0\) on \(0<r<1\),
\(\phi(1)=0\) and \(\phi'(1)<0\), we have \(\phi\Delta\phi>0\) in
an interior collar of the wall. Take a nonzero smooth radial cutoff
\(\chi\ge0\) supported in that collar and put
\[
 a(r)=1+\frac{\int_D|\nabla\phi|^2\,dx}
                  {\int_D\chi\phi\Delta\phi\,dx}\,\chi(r).
\]
The zero trace of \(\phi\) gives \(\int_D\nabla\phi\cdot\nabla(a\phi)=0\)
by integration by parts. The spatial profile has a single time factor,
so the same weight gives \(\nabla F(t)\perp_{\BJ}\nabla(aF(t))\)
for every \(t\ge0\).

\smallskip
\noindent\textbf{Subcritical reverse estimates for the full velocity.}
We now give a sufficient condition for continuation in three dimensions
without a symmetry hypothesis. All velocity and gradient norms below
use the Euclidean and Frobenius norms, respectively.

\begin{proposition}\label{prop:velocity-continuation}
Let \(u\in C([0,T);H^3(\R^3))\) be a smooth solution of the unforced
incompressible Navier--Stokes equations on \([0,T)\times\R^3\), with
viscosity \(\nu>0\) and \(T<\infty\).
If there exist \(0\le t_0<T\), \(2<q<6\) and \(0<C<\infty\)
such that, for almost every \(t\in(t_0,T)\),
\begin{equation}\label{eq:velocity-subcritical}
 \displaystyle
 \left(\int_{\R^3}|\nabla u(t,x)|^2\,dx\right)^{1/2}
 \le C\left(\int_{\R^3}|u(t,x)|^q\,dx\right)^{1/q},
\end{equation}
then \(u\) extends smoothly past \(T\).
The same continuation conclusion holds under the
following fixed-weight condition. For one \(a\in W^{1,\infty}(\R^3)\)
with \(m=\inf a>0\), assume
\[
 \int_{\R^3}\sum_{i,j=1}^3\partial_j u_i(t,x)\,
                    \partial_j(a(x)u_i(t,x))\,dx\le0
\]
for almost every \(t\in(t_0,T)\). Equality is the Birkhoff--James
condition for the two matrix gradients in \(L^2(\R^3;\R^{3\times3})\).
\end{proposition}

\begin{proof}
For \(2<q<6\), put \(\theta=3(1/2-1/q)\in(0,1)\).
Interpolation and
\(\norm{u}_6\le S\norm{\nabla u}_2\) give
\[
 \norm{\nabla u(t)}_2\le CS^\theta
       \norm{u(t)}_2^{1-\theta}\norm{\nabla u(t)}_2^\theta,\qquad
 \norm{\nabla u(t)}_2\le (CS^\theta)^{1/(1-\theta)}\norm{u(0)}_2,
\]
where the second estimate uses the kinetic energy
bound in \eqref{eq:ns-kinetic}. This bounds
the \(H^1\) norm up to \(T\), so the whole-space \(H^1\) continuation
theory applies \cite[Corollary~5.8]{Tao2013}.
For the weight criterion, expand the product and use Cauchy--Schwarz
to obtain \(m\norm{\nabla u(t)}_2^2\le
\norm{\nabla a}_\infty\norm{u(t)}_2\norm{\nabla u(t)}_2\).
This is a direct bound at the lower endpoint \(q=2\):
\[
 \|\nabla u(t)\|_2
 \le \frac{\|\nabla a\|_\infty}{m}\|u(t)\|_2
 \le \frac{\|\nabla a\|_\infty}{m}\|u(0)\|_2 .
\]
The resulting uniform \(H^1\) bound gives continuation.
\end{proof}

Under the same one-sided pairing condition, if
\(\nabla a\in L^{2q/(q-2)}(\R^3)\) for some \(2<q<6\), then
H\"older's inequality also gives the estimate in
\eqref{eq:velocity-subcritical} directly with
\(C=m^{-1}\|\nabla a\|_{L^{2q/(q-2)}(\R^3)}\).
In particular, a fixed weight whose gradient is bounded and compactly
supported yields the reverse estimate throughout \(2<q<6\).

For this three-dimensional reverse quotient, \(q=6\)
is the Sobolev scaling-critical exponent. At \(q=6\) the kinetic energy
drops out of the interpolation estimate.
For \(u_\lambda(x)=\lambda u(\lambda x)\), both
\(\norm{\nabla u_\lambda}_2\) and \(\norm{u_\lambda}_6\) grow as
\(\lambda^{1/2}\). Thus bounding their ratio leaves this
concentration of spatial profiles possible. An implication from the
critical estimate alone to continuation would require a further
argument using the evolution.

\smallskip
\noindent\textbf{Axisymmetric flow without swirl.}
When \(u^\theta=0\), set
\(\eta=\omega^\theta/r=(\partial_z u^r-\partial_r u^z)/r\).
The vorticity equation \cite[Eq.~(1.8)]{GallaySverak2015} becomes
\[
 \partial_t\eta+u^r\partial_r\eta+u^z\partial_z\eta
 =\nu\bigl(\partial_r^2+3r^{-1}\partial_r+\partial_z^2\bigr)\eta .
\]
The diffusion on the right is the Laplacian on \(\R^4\times\R\)
for functions radial in the first four variables. The maximum
principle preserves the sign of one-signed initial relative vorticity.
For smooth decaying solutions, multiplication by \(\eta\) and
integration against \(dx=2\pi r\,dr\,dz\) give
\begin{equation}\label{eq:no-swirl-energy}
 \frac12\frac{d}{dt}\int_{\R^3}\eta^2\,dx
 +\nu\int_{\R^3}|\nabla\eta|^2\,dx
 +2\pi\nu\int_{\R}|\eta(t,0,z)|^2\,dz=0 .
\end{equation}
In fact, the transport term vanishes by incompressibility and
the extra radial diffusion contributes
\(4\pi\nu\int\eta\,\partial_r\eta\,dr\,dz
=-2\pi\nu\int\eta(t,0,z)^2\,dz\).
Proposition~\ref{prop:profile-weight} now shows that, for one-signed
initial \(\eta\), every nonzero \(C^2\cap H^2(\R^3)\) time slice
admits a smooth positive axisymmetric weight with
\(\nabla\eta(t)\perp_{\BJ}\nabla(a_t\eta(t))\).
The classical global theory also covers sign-changing data
\cite{GallaySverak2015}.

\smallskip
\noindent\textbf{Two planar vorticity examples.}
For \(D=\R^2\) or \(\mathbb T^2\), taking the curl of
\eqref{eq:ns-system} gives the scalar vorticity equation
\cite{GallagherGallay2005}
\begin{equation}\label{eq:ns-planar-vorticity}
 \partial_t\omega+u_1\partial_1\omega+u_2\partial_2\omega
   =\nu(\partial_1^2+\partial_2^2)\omega,\qquad
 \omega=\partial_1u_2-\partial_2u_1 .
\end{equation}
For a smooth decaying or periodic solution, incompressibility gives
\(\int_D(u\cdot\nabla\omega)\omega\,dx=0\). Hence
\[
 \frac12\frac{d}{dt}\int_D\omega^2\,dx
  +\nu\int_D|\nabla\omega|^2\,dx=0 .
\]

The planar argument uses a finite exponent \(2<q<\infty\). Theorems
\ref{thm:relaxation} and \ref{thm:whole-space-relaxation} use the
critical Lebesgue exponent \(2n/(n-2)\) in dimension \(n\ge3\).
In dimension two the first-order \(L^2\) Sobolev scale has no finite
critical Lebesgue exponent. For the planar estimate take \(D=\R^2\) or \(\mathbb T^2\),
\(\omega\in H^1(D)\cap L^q(D)\), and a bounded positive
Lipschitz weight with \(m=\inf_D a>0\) and
\(\nabla a\in L^{2q/(q-2)}(D)\). If
\(\int_D\nabla\omega\cdot\nabla(a\omega)\,dx\le0\), expansion of the
product gives
\[
 m\|\nabla\omega\|_{L^2(D)}^2
 \le \|\nabla\omega\|_{L^2(D)}\|\omega\nabla a\|_{L^2(D)}.
\]
Since \(1/2=1/q+(q-2)/(2q)\), H\"older's inequality yields
\begin{equation}\label{eq:planar-reverse}
 \displaystyle
 \|\nabla\omega\|_{L^2(D)}
 \le \frac{\|\nabla a\|_{L^{2q/(q-2)}(D)}}{\inf_D a}
       \|\omega\|_{L^q(D)},\qquad 2<q<\infty.
\end{equation}
Exact Birkhoff--James orthogonality gives equality to zero in the
pairing hypothesis. At the lower H\"older endpoint \(q=2\), the same
argument uses \(\|\nabla a\|_{L^\infty(D)}\). The two flows below give
explicit tests of when the planar pairing vanishes.

First let \(s=1+4\nu t\), \(y=|x|^2/s\), and take
\[
 \psi(t,x)=s^{-1}e^{-y},\qquad
 u=(-\partial_2\psi,\partial_1\psi),\qquad \omega=\Delta\psi .
\]
The stream function solves the heat equation, and the tangential
velocity transports the radial vorticity along its level sets.
This gives a smooth finite-energy Navier--Stokes solution on \(\R^2\).
Direct differentiation gives
\[
 \omega=4s^{-2}(y-1)e^{-y},\qquad
 \Delta\omega=16s^{-3}(y^2-4y+2)e^{-y}.
\]
Their product is positive on \(y>2+\sqrt2\). Choose
\(0\le\chi\in C_c^\infty(\{|x|^2>2+\sqrt2\})\) nonzero and define
\(a_0\) by \eqref{eq:profile-cutoff} with \(f=\omega(0)\) and
cutoff \(\chi\). Then for every \(t\ge0\),
\begin{equation}\label{eq:planar-weight}
 a_t(x)=a_0(x/\sqrt s),\qquad
 \nabla\omega(t)\perp_{\BJ}\nabla(a_t\omega(t)),\qquad
 \norm{\nabla a_t}_\infty=s^{-1/2}\norm{\nabla a_0}_\infty .
\end{equation}
Indeed \(\omega(t,x)=s^{-2}\omega(0,x/\sqrt s)\), so changing
variables multiplies the initial gradient pairing by \(s^{-4}\).
The weights \(a_t\) have compactly supported gradients. For each
finite \(q>2\), the coefficient in \eqref{eq:planar-reverse} equals
\[
 \frac{\|\nabla a_t\|_{L^{2q/(q-2)}(\R^2)}}{\inf_{\R^2}a_t}
 =s^{-1/q}\frac{\|\nabla a_0\|_{L^{2q/(q-2)}(\R^2)}}{\inf_{\R^2}a_0}.
\]

On the \(2\pi\)-periodic torus the globally smooth shear
\(u(t,x)=(e^{-\nu t}\sin x_2,0)\) has
\(\omega=-e^{-\nu t}\cos x_2\) and \(-\Delta\omega=\omega\).
Periodic integration by parts gives both the lower-order norm identity
and the strict sign of the pairing for every smooth positive periodic
weight \(a\) at every finite time,
\begin{equation}\label{eq:planar-shear}
 \int_{\mathbb T^2}\nabla\omega\cdot\nabla(a\omega)\,dx
   =\int_{\mathbb T^2}a\omega^2\,dx>0,\qquad
 \norm{\nabla\omega}_{L^2(\mathbb T^2)}=\norm{\omega}_{L^2(\mathbb T^2)}.
\end{equation}


\begin{thebibliography}{18}

\bibitem{James1947}
James, R.C.: Orthogonality and linear functionals in normed linear
spaces. Trans. Amer. Math. Soc. \textbf{61}, 265--292 (1947).
\href{https://doi.org/10.1090/S0002-9947-1947-0021241-4}{DOI}

\bibitem{LangMihula2023}
Lang, J., Mihula, Z.: Different degrees of non-compactness for optimal
Sobolev embeddings. J. Funct. Anal. \textbf{284}, 109880 (2023).
\href{https://doi.org/10.1016/j.jfa.2023.109880}{DOI}

\bibitem{Talenti1976}
Talenti, G.: Best constant in Sobolev inequality. Ann. Mat. Pura Appl.
\textbf{110}, 353--372 (1976).
\href{https://doi.org/10.1007/BF02418013}{DOI}

\bibitem{MuruganSaloffCoste2017}
Murugan, M., Saloff-Coste, L.: Davies' method for anomalous diffusions.
Proc. Amer. Math. Soc. \textbf{145}, 1793--1804 (2017).
\href{https://doi.org/10.1090/proc/13324}{DOI}

\bibitem{Vazquez2006}
V\'azquez, J.L.: Smoothing and Decay Estimates for Nonlinear Diffusion
Equations. Equations of Porous Medium Type. Oxford Lecture Series in
Mathematics and Its Applications, vol.~33. Oxford University Press,
Oxford (2006).
\href{https://doi.org/10.1093/acprof:oso/9780199202973.001.0001}{DOI}

\bibitem{Palasek2026}
Palasek, S.: A logarithmic regularity criterion for the axisymmetric
Navier--Stokes equations. Preprint (2026).
\href{https://arxiv.org/abs/2610.29012}{arXiv:2610.29012}

\bibitem{ChenFangZhang2017}
Chen, H., Fang, D., Zhang, T.: Regularity of 3D axisymmetric
Navier--Stokes equations. Discrete Contin. Dyn. Syst. \textbf{37},
1923--1939 (2017).
\href{https://doi.org/10.3934/dcds.2017081}{DOI}

\bibitem{Tao2013}
Tao, T.: Localisation and compactness properties of the
Navier--Stokes global regularity problem. Anal. PDE \textbf{6},
25--107 (2013).
\href{https://doi.org/10.2140/apde.2013.6.25}{DOI}

\bibitem{GallaySverak2015}
Gallay, T., \v{S}ver\'ak, V.: Remarks on the Cauchy problem for the
axisymmetric Navier--Stokes equations. Confluentes Math. \textbf{7}(2),
67--92 (2015).
\href{https://doi.org/10.5802/cml.25}{DOI}

\bibitem{Leray1934}
Leray, J.: Sur le mouvement d'un liquide visqueux
emplissant l'espace. Acta Math. \textbf{63}, 193--248 (1934).
\href{https://doi.org/10.1007/BF02547354}{DOI}

\bibitem{Serrin1962}
Serrin, J.: On the interior regularity of weak
solutions of the Navier--Stokes equations. Arch. Ration. Mech. Anal.
\textbf{9}, 187--195 (1962).
\href{https://doi.org/10.1007/BF00253344}{DOI}

\bibitem{EscauriazaSereginSverak2003}
Escauriaza, L., Seregin, G.A., \v{S}ver\'ak, V.:
\(L_{3,\infty}\)-solutions of the Navier--Stokes equations and
backward uniqueness. Russian Math. Surveys \textbf{58}, 211--250 (2003).
\href{https://doi.org/10.1070/RM2003v058n02ABEH000609}{DOI}

\bibitem{KochNadirashviliSereginSverak2009}
Koch, G., Nadirashvili, N., Seregin, G.A.,
\v{S}ver\'ak, V.: Liouville theorems for the Navier--Stokes equations
and applications. Acta Math. \textbf{203}, 83--105 (2009).
\href{https://arxiv.org/abs/0709.3599}{arXiv:0709.3599}

\bibitem{GallagherGallay2005}
Gallagher, I., Gallay, T.: Uniqueness for the
two-dimensional Navier--Stokes equation with a measure as initial
vorticity. Math. Ann. \textbf{332}, 287--327 (2005).
\href{https://doi.org/10.1007/s00208-004-0627-x}{DOI}

\bibitem{GladialiGrossiNeves2013}
Gladiali, F., Grossi, M., Neves, S.L.N.: Nonradial solutions for the
H\'enon equation in \(\R^N\). Adv. Math. \textbf{249}, 1--36 (2013).
\href{https://doi.org/10.1016/j.aim.2013.07.022}{DOI}

\bibitem{KaminVazquez1988}
Kamin, S., V\'azquez, J.L.: Fundamental solutions and asymptotic
behaviour for the \(p\)-Laplacian equation. Rev. Mat. Iberoam.
\textbf{4}, 339--354 (1988).
\href{https://doi.org/10.4171/RMI/77}{DOI}

\bibitem{ZhangZhang2014}
Zhang, P., Zhang, T.: Global axisymmetric solutions to three-dimensional
Navier--Stokes system. Int. Math. Res. Not. IMRN \textbf{2014},
610--642 (2014).
\href{https://doi.org/10.1093/imrn/rns232}{DOI}

\bibitem{Konig2023}
K\"onig, T.: On the sharp constant in the Bianchi--Egnell stability
inequality. Bull. Lond. Math. Soc. \textbf{55}, 2070--2075 (2023).
\href{https://doi.org/10.1112/blms.12837}{DOI}

\end{thebibliography}
\end{document}